\documentclass[11pt]{article}
\usepackage{}
\usepackage{amssymb}
\usepackage{epsfig}
\def\be{\begin{equation}}
\def\ee{\end{equation}}
\def\bea{\begin{eqnarray}}
\def\eea{\end{eqnarray}}
\def\bes{\begin{eqnarray*}}
\def\ees{\end{eqnarray*}}

\def\nn{\nonumber}
\def\<{\langle}
\def\>{\rangle}
\def\lb{\label}

\def\R{{\bf R}}

\def\Z{{\bf Z}}

\def\Sg{{\Sigma}}

\def\hb{\vrule height0.18cm width0.14cm $\,$}

\title{ On the number of P-invariant closed characteristics on partially symmetric compact convex hypersurfaces
in ${\bf R}^{2n}$}
\author{Hui
Liu$^{1}$, \thanks{Partially supported by China Postdoctoral Science
Foundation No.2013M540512.
E-mail: huiliu@ustc.edu.cn} \qquad  Duanzhi Zhang$^{2}$, \thanks{Partially supported by NSF of China (10801078, 11171341, 11271200) and LPMC of Nankai
University. E-mail: zhangdz@nankai.edu.cn}\\ \\
$^{1}$ School of Mathematical Sciences, University of Science and Technology of China, \\Hefei, Anhui 230026, P. R. China\\
$^{2}$ School of Mathematical Sciences and LPMC, Nankai University,\\ Tianjin 300071, P. R. China }
\date{}

\begin{document}

\maketitle

\begin{abstract}
{\it  In this paper, let $n\geq2$ be an integer,
$P=diag(-I_{n-\kappa},I_\kappa,-I_{n-\kappa},I_\kappa)$ for some
integer $\kappa\in[0, n)$, and $\Sigma \subset {\bf R}^{2n}$ be a
partially symmetric compact convex hypersurface, i.e., $x\in \Sigma$ implies $Px\in\Sigma$.
We prove that if $\Sigma$ is $(r,R)$-pinched with
$\frac{R}{r}<\sqrt{2}$, then there exist at least $n-\kappa$ geometrically distinct P-symmetric closed
characteristics on $\Sigma$, as a consequence, $\Sigma$ carry at least $n$ geometrically distinct P-invariant closed
characteristics. }
\end{abstract}



\renewcommand{\theequation}{\thesection.\arabic{equation}}
\renewcommand{\thefigure}{\thesection.\arabic{figure}}

\setcounter{equation}{0}%\setcounter{figure}{0}
\section{Introduction and main results}%{Section 1}
Let $\Sigma$ be a $C^2$ compact hypersurface in ${\bf R}^{2n}$,
bounding a strictly convex compact set $U$ with non-empty interior, where $n\geq2$.
We denote the set of all such hypersurfaces by $\mathcal{H}(2n)$.
Without loss of generality, we suppose $U$ contains the origin. We
consider closed characteristics $(\tau,y)$ on $\Sigma$, which are
solutions of the following problem \be
\left\{\matrix{\dot{y}(t)=JN_\Sigma(y(t)),~y(t)\in \Sg,~\forall~t\in
{\bf R}, \cr
               y(\tau)=y(0), ~~~~~~~~~~~~~~~~~~~~~~~~~~~~~~~~~\cr }\right. \lb{1.1}\ee
where $J=\left(
              \begin{array}{cc}
                0 & -I_{n} \\
                I_{n} & 0 \\
              \end{array}
            \right)
$, $I_{n}$ is the identity matrix in ${\bf R}^n$ and
$\mathit{N}_\Sigma(y)$ is the outward normal unit vector  of
$\Sigma$ at $y$ normalized by the condition
$\mathit{N}_\Sigma(y)\cdot y=1$. Here $a\cdot b$ denotes the
standard inner product of $a, b \in {\bf R}^{2n}$. A closed
characteristic $(\tau, y)$ is {\it prime} if $\tau$ is the minimal
period of $y$. Two closed characteristics $(\tau,x)$ and
$(\sigma,y)$ are {\it geometrically distinct}, if $x({\bf
R})\not=y({\bf R})$.  We denote by $\mathcal{T}(\Sg)$ the set of all
geometrically distinct closed characteristics on $\Sg$.

There is a long standing conjecture on the number of closed
characteristics on compact convex hypersurfaces in ${\bf
R}^{2n}$\bea \,^\#\mathcal{T}(\Sg)\geq n,~~~\forall
~\Sigma\in\mathcal{H}(2n)\lb{1.2}\eea Since the pioneering works \cite{Rab1}
of P. Rabinowitz and \cite{Wei1} of A. Weinstein in 1978 on the
existence of at least one closed characteristic on every
hypersurface in $\mathcal {H}(2n)$, the existence of multiple closed
characteristics on $\Sigma\in\mathcal{H}(2n)$ has been deeply
studied by many mathematicians. In particular, I. Ekeland and J. Lasry
proved (\ref{1.2}) under a $\sqrt{2}$-pinched condition(cf. \cite{EkL1}), and
under a $\sqrt{3}$-pinched condition(cf. \cite{Gir1}), M. Girardi proved that
there exist $n$ geometrically distinct symmetric closed
characteristics on $\Sigma \in\mathcal{H}(2n)$ which is symmetric with respect to the origin,
i.e., $x\in\Sigma$ implies $-x\in \Sigma$. In 1987-1988 I. Ekeland-L.
Lassoued, I. Ekeland-H. Hofer, and A. Szulkin (cf. \cite{EkL2},
\cite{EkH1}, \cite{Szu1}) proved \bea \,^\#\mathcal{T}(\Sg)\geq
2,~~~\forall ~\Sigma\in\mathcal{H}(2n)\lb{1.3}\eea In \cite{LoZ1} of 2002,
Y. Long and C. Zhu proved \bea \,^\#\mathcal{T}(\Sg)\geq
[\frac{n}{2}]+1,~~~\forall~ \Sigma\in\mathcal{H}(2n)\lb{1.4}\eea  In
\cite{WHL1}, the authors proved the conjecture for $n = 3$. In
\cite{LLZ1}, the the authors proved the conjecture when
$\Sigma\in\mathcal{H}(2n)$ is symmetric with respect to the origin,
i.e., $x\in\Sigma$ implies $-x\in \Sigma$. In \cite{DoL1} of 2004,
Y. Dong and Y. Long studied the multiplicity of closed
characteristics on partially symmetric convex hypersurfaces in ${\bf
R}^{2n}$.

For any $s_i, t_i\in{\bf R}^{k_i}$ with $i=1,2$, we denote by $(s_1,
t_1)\diamond(s_2, t_2)=(s_1, s_2, t_1, t_2)$. Fixing an integer
$\kappa$ with $0\leq \kappa < n$, let
$P=diag(-I_{n-\kappa},I_\kappa,-I_{n-\kappa},I_\kappa)$ and
$\mathcal {H}_\kappa(2n)=\{\Sigma\in\mathcal {H}(2n)\mid x\in\Sigma
~implies~ Px\in\Sigma\}$. For $\Sigma\in\mathcal {H}_\kappa(2n)$,
let $\Sigma(\kappa)=\{z\in {\bf R}^{2\kappa}\mid 0\diamond
z\in\Sigma\}$, where 0 is the origin in ${\bf R}^{2n-2\kappa}$. As
in \cite{DoL1}, A closed characteristic $(\tau, y)$ on $\Sigma \in
\mathcal {H}_\kappa(2n)$ is {\it P-asymmetric} if $y({\bf R})\cap
Py({\bf R})=\emptyset$, it is {\it P-symmetric} if $y({\bf
R})=Py({\bf R})$ with $y=y_1\diamond y_2$ and $y_1\neq 0$, or it is
{\it P-fixed} if $y({\bf
R})=Py({\bf R})$ and $y=0\diamond y_2$, where
$y_1\in{\bf R}^{2(n-\kappa)}$, $y_2\in{\bf R}^{2\kappa}$. We call a
closed characteristic $(\tau, y)$ is {\it P-invariant} if $y({\bf
R})=Py({\bf R})$. Then a P-invariant
closed characteristic is P-symmetric or
P-fixed. Note that if we set $\kappa=0$, then the P-symmetric closed
characteristic is just symmetric and the P-fixed closed
characteristics vanish.

It is very interesting to consider closed characteristics on hypersurfaces with
special symmetries. Recently, W. Wang proved that there exist two
symmetric closed characteristics on any compact convex hypersurface which is symmetric with respect to the origin (cf. \cite{Wan1}),
the first author of this paper proved that there exist at least two
geometrically distinct P-invariant closed
characteristics on any $\Sigma\in\mathcal {H}_\kappa(2n)$ (cf. \cite{Liu1}) and the second author of this paper
proved that there exist at least two
geometrically distinct P-cyclic closed
characteristics on P-cyclic compact convex hypersurfaces (cf. \cite{Zha1}).
Thus whether $\Sigma\in\mathcal {H}_\kappa(2n)$ carries more P-invariant closed
characteristics becomes an interesting problem. Motivated by \cite{EkL1}, \cite{Liu1} and \cite{Zha1},
we prove the following results in this paper.

{\bf Theorem 1.1.} {\it Assume $\Sigma\in\mathcal {H}_\kappa(2n)$ and $0<r \leq |x|\leq R,~\forall~x \in \Sigma$
with $\frac{R}{r}<\sqrt{2}$. Then there exist at least $n-\kappa$ geometrically distinct P-symmetric closed
characteristics $(\tau_i, y_i)$ on $\Sigma$, where $\tau_i$ is the minimal period of $y_i$,
and the actions $A(\tau_i, y_i)$ satisfy:  $\pi r^2\leq A(\tau_i, y_i)\leq \pi R^2, \forall 1\leq i\leq n-\kappa$.}

Combining Theorem 1.1 and the result proved by I. Ekeland and J. Lasry in 1980, we have a direct consequence:

{\bf Theorem 1.2.} {\it Assume $\Sigma\in\mathcal {H}_\kappa(2n)$ and $0<r \leq |x|\leq R,~\forall~x \in \Sigma$
with $\frac{R}{r}<\sqrt{2}$. Then $\Sigma$ carries at least $n$ geometrically distinct P-invariant closed
characteristics $(\tau_i, y_i)$, where $\tau_i$ is the minimal period of $y_i$,
and the actions $A(\tau_i, y_i)$ satisfy:  $\pi r^2\leq A(\tau_i, y_i)\leq \pi R^2, \forall 1\leq i\leq n$..}

Here the action of a closed characteristic
$(\tau,y)$ is defined by (cf. P.190 of \cite{Eke1})\[
A(\tau,y)=\frac{1}{2}\int_0^{\tau}{(Jy\cdot \dot{y})}dt.\]

In this paper, let ${\bf N}, {\bf N}_0, {\bf Z}, {\bf Q}, {\bf R}$
and ${\bf C}$ denote the sets of natural integers, non-negative
integers, integers, rational numbers, real numbers and complex
numbers respectively. Denote by $a\cdot b$ and $|a|$ the standard
inner product and norm in ${\bf R}^{2n}$. Denote by $\langle
\cdot,\cdot \rangle$ and $\|\cdot\|$ the standard $L^2$-inner
product and $L^2$-norm. For a set $A$, we denote by $\,^\#A$ the
number of elements in $A$. For an $S^1$-space $X$, we denote by
$X_{S^1}$ the homotopy quotient of $X$ by $S^1$, i.e.,
$X_{S^1}=S^\infty\times_{S^1}X$, where $S^\infty$ is the unit sphere
in an infinite dimensional {\it complex} Hilbert space.

\setcounter{equation}{0}
\section{A variational structure for P-invariant closed characteristics}
In the rest of this paper, we fix a $\Sigma\in
\mathcal{H}_\kappa(2n)$. Note that a prime closed characteristic $(\tau, y)$ is P-symmetric if and only if it satisfies
the problem
\be  \left\{\matrix{
      \dot{y}(t)=JN_{\Sigma}(y(t)),  y(t)\in \Sg,~~~~~~~~~~~~~~~~ \cr
              y(\frac{\tau}{2})=P y(0), y(t)=y_1(t)\diamond y_2(t), y_1\neq 0,       \cr}\right.  \lb{2.1}\ee
and a prime closed characteristic $(\tau, y)$ is P-fixed if and only if it satisfies
the problem
\be  \left\{\matrix{
      \dot{y}(t)=JN_{\Sigma}(y(t)),  y(t)\in \Sg, \cr
              y(\tau)= y(0), y(t)=0\diamond y_2(t),       \cr}\right.  \lb{2.2}\ee
where $y_1(t)\in{\bf R}^{2(n-\kappa)}, y_2(t)\in{\bf R}^{2\kappa}, \forall\;t\in {\bf R}$. Here we also note that
even iterate $(2m\tau,y)$ of any prime P-symmetric closed characteristic $(\tau, y)$ does not satisfy the equation (\ref{2.1}).

In this section, we transform the problems (\ref{2.1})-(\ref{2.2}) into a fixed
period problem of a Hamiltonian system and then study its
variational structure.

As on Page 199 of \cite{Eke1}, we choose some $\alpha\in (1,2)$ and associate with $U$ a convex function
$H(x)=j(x)^\alpha, \forall x\in\R^{2n}$, where $j: {\bf R}^{2n}\rightarrow {\bf R}$ is the gauge function of
$\Sigma$, that is, $j(\lambda x)=\lambda$ for $x\in\Sigma$ and
$\lambda\geq 0$. Then we have the following:

{\bf Proposition 2.1.} {\it Consider the following problem
\bea\left\{\matrix{ \dot{x}(t)=JH^{\prime}(x(t)), \cr
x(\frac{1}{2})=P x(0).~~~~\cr }\right.\lb{2.3}\eea
Then solutions of (\ref{2.3}) are $x\equiv 0$ and $x={(\frac{\tau}{\alpha})}^{-\frac{1}{2-\alpha}} y(\tau
t)$, where $(\tau,y)$ is a solution of (\ref{2.1}), or
$(\frac{\tau}{2},y)$ is a P-fixed closed characteristic. In
particular, nonzero solutions of (\ref{2.3}) are in one to one
correspondence with P-invariant closed characteristics.}

{\it Proof.} Clearly $x\equiv 0$ is the unique constant solution of (\ref{2.3}).
Suppose $x(t)$ is a nonconstant solution of (\ref{2.3}), then
$H(x(t))=(j(x(t)))^\alpha=const$. Let $\rho=j(x(t))$ and
$y(t)=\rho^{-1}x\left(\frac{\rho^{2-\alpha}}{\alpha}t\right)$.
Then $j(y)=\rho^{-1}j(x)=\rho^{-1}\rho=1$, hence $y(\R)\subset \Sigma$.
Moreover, noticing the fact that $j^\prime(y)=N_\Sigma(y)$ and $j(\lambda x)=\lambda j(x)$ for all
$x\in \R^{2n}\setminus \{0\}$ and $\lambda\in \R^+$, we have $\dot y(t)=JN_\Sigma(y(t))$ and
$y(\frac{\alpha}{2\rho^{2-\alpha}})=P y(0)$ by (\ref{2.3}).
Let $\tau=\frac{\alpha}{\rho^{2-\alpha}}$, then $(\tau, y)$ is a closed characteristic satisfying (\ref{2.1}), or
$(\frac{\tau}{2}, y)$ is a P-fixed closed characteristic. The other side of the proposition can be proved similarly
and thus is omitted.      \hfill\hb

In the following, we use the Clarke-Ekeland dual action principle to problem (\ref{2.3}). As usual, let
$G$ be the Fenchel transform of $H$ defined by $G(y)=\sup{\{x\cdot y-H(x) \mid x\in{\bf R}^{2n}\}}$.
Then $G$ is strictly convex.
As in Section 3 of \cite{DoL2} (also Section 2 of \cite{Liu1}), let \bea L^2_\kappa\left(0,
\frac{1}{2}\right)&=&\{u=u_1\diamond u_2 \in L^2((0,\frac{1}{2}),
{\bf R}^{2n})\mid u_1\in L^2((0,\frac{1}{2}), {\bf
R}^{2n-2\kappa}),\nn\\&& u_2\in L^2((0,\frac{1}{2}), {\bf
R}^{2\kappa}), u(\frac{1}{2})=Pu(0),
\int_0^{\frac{1}{2}}u_2(t)dt=0\}.\lb{2.4}\eea
Define a linear operator $\Pi_\kappa: L^2_\kappa\left(0,
\frac{1}{2}\right)\rightarrow L^2_\kappa\left(0, \frac{1}{2}\right)$
by\bea (\Pi_\kappa u)(t)&=&x_1(t)\diamond x_2(t),\lb{2.5}\\x_1(t)&=&\int_0^t
u_1(\tau)d\tau-\frac{1}{2}\int_0^{\frac{1}{2}}u_1(\tau)d\tau,\lb{2.6}\\x_2(t)&=&\int_0^t
u_2(\tau)d\tau-2\int_0^{\frac{1}{2}}dt\int_0^t u_2(\tau)d\tau,\lb{2.7}\eea
for any $u=u_1\diamond u_2\in L^2_\kappa\left(0,
\frac{1}{2}\right)$. Then by Lemma 2.4 of \cite{Liu1}, we know $\Pi_\kappa$ is a compact operator from
$L^2_\kappa\left(0, \frac{1}{2}\right)$ into itself and
$J\Pi_\kappa$ is self-adjoint. Now the dual action functional on $L^2_\kappa\left(0,
\frac{1}{2}\right)$ is defined by
\be\Psi(u)=\int_0^{\frac{1}{2}}\left(\frac{1}{2}J u \cdot
\Pi_\kappa u+G(-Ju)\right)dt.\lb{2.8}\ee
Since $H$ is $\alpha$-homogeneous, we have $H(x)\leq r^{-\alpha} |x|^\alpha$ for some positive constant,
then by the definition of $G$, we get $G(y)\geq \frac{1}{\beta}(\frac{r^\alpha}{\alpha})^{\beta-1}|y|^\beta$,
where $\beta>2$ satisfying $\alpha^{-1}+\beta^{-1}=1$, cf. also (7) and (22) in Section V.2 of \cite{Eke1}.
Then by the same proofs of Propositions 2.5 and 2.6 in \cite{Liu1}, we have

{\bf Lemma 2.2.} {\it  The functional $\Psi$ is $C^{1,1}$ and is bounded
from below on $L^2_\kappa\left(0, \frac{1}{2}\right)$ and satisfies the
Palais-Smale condition. Suppose $x$ is a solution of (\ref{2.3}).
Then $u =\dot{x}$ is a critical point of $\Psi$. Conversely,
suppose $u$ is a critical point of $\Psi$. Then there exists a
unique $\xi\in{\bf R}^{2\kappa}$ such that $\Pi_\kappa u-0\diamond\xi$ is a
solution of (\ref{2.3}). In particular, solutions of (\ref{2.3}) are
in one to one correspondence with critical points of $\Psi$.}

Note that we can identify $L^2_\kappa\left(0,
\frac{1}{2}\right)$ with the space $\{u\in L^2(\R/\Z,
{\bf R}^{2n})\mid u|_{(0,1/2)}\in L^2_\kappa(0,
\frac{1}{2}), u(t+\frac{1}{2})=Pu(t) \}$. Then we have a natural $S^{1}$-action on $L^2_\kappa\left(0,
\frac{1}{2}\right)$ defined by $\theta*u(t) = u(\theta+t)$, for all
$\theta \in S^{1}\equiv {\bf R}/{\bf Z}$ and $t \in {\bf R}$. Geometrically, under a partial reflection $P$,
$u([0, 1/2])$ equals to $u([\theta, 1/2+\theta])$ for all $\theta \in S^{1}$ and $\Psi$ is invariant under the
partial reflection $P$, i.e., $\Psi(Pu)=\Psi(u)$. As Lemma 2.8 of \cite{Liu1}, we have

{\bf Lemma 2.3.} {\it The functional $\Psi$ is
$S^{1}$-invariant.}

{\bf Definition 2.4.} {\it Define the $S^1$-orbit of a point $u\in L^2_\kappa\left(0,
\frac{1}{2}\right)$ by\bea \mathcal {O}
(u):=\{\theta*u\mid \theta\in S^1\}.\nn\eea
We shall say that the $S^1$-action is free at $u$ if the map $\theta\rightarrow \theta*u$ is injective; that is, a
homeomorphism of $S^1$ onto $\mathcal {O}
(u)$. If $\Omega\subset L^2_\kappa\left(0,
\frac{1}{2}\right)$ is an invariant subset, the $S^1$-action is free in $\Omega$ if it is free at every point $u\in\Omega$.}

In our case, the $S^1$-action will be free at $u$ if and only if $u$ has minimal period 1. Note that
by Lemma 2.3, the level sets $\Psi^\gamma:=\{u\in L^2_\kappa\left(0,
\frac{1}{2}\right)\mid \Psi(u)<\gamma\}$ are $S^1$-invariant.

Now by Proposition 2.1 and Lemma 2.2, we assume $u=\dot{x}$ is critical point of $\Psi$ corresponding to a solution $(\tau, y)$ of (\ref{2.1}), or
a P-fixed closed characteristic $(\frac{\tau}{2},y)$, then $x=\Pi_\kappa u-0\diamond\xi$ for some $\xi\in{\bf R}^{2\kappa}$,
the corresponding critical value of $\Psi$ is given by:

{\bf Lemma 2.5.} {\it $\Psi(u)=-\frac{1}{2}(1-\frac{\alpha}{2})(\frac{\tau}{\alpha})^{-\frac{\alpha}{2-\alpha}}$.}

{\bf Proof.} Since $-J u=-J \dot{x}=H^{\prime}(x)$, then $G(-J u)=x\cdot H^{\prime}(x)-H(x)$. Thus by direct computations,
we have \bea
\Psi(u)&=&\int_0^{\frac{1}{2}}\left(\frac{1}{2}J \dot{x} \cdot
(x+0\diamond\xi)+G(-J \dot{x})\right)dt\nn\\&=&\int_0^{\frac{1}{2}}\left(\frac{1}{2}J \dot{x} \cdot
x+x\cdot H^{\prime}(x)-H(x)\right)dt\nn\\&=&\int_0^{\frac{1}{2}}\left(\frac{1}{2}x\cdot H^{\prime}(x)-H(x)\right)dt.\lb{2.9}\eea
Note that $H(\lambda x)=\lambda^\alpha H(x)$ for all $\lambda\geq 0$, then $x\cdot H^{\prime}(x)=\alpha H(x)$ and (\ref{2.9})
becomes  \bea
\Psi(u)&=&\int_0^{\frac{1}{2}}\left((\frac{1}{2}\alpha-1)H(x)\right)dt\nn\\
&=&-\frac{1}{2}(1-\frac{\alpha}{2})(\frac{\tau}{\alpha})^{-\frac{\alpha}{2-\alpha}},\nn\eea
here we used the fact that $x={(\frac{\tau}{\alpha})}^{-\frac{1}{2-\alpha}} y(\tau
t)$.\hfill\hb

\setcounter{equation}{0}
\section{ Proofs of the main theorems}
In this section, we assume $0<r \leq |x|\leq R,~\forall~x \in \Sigma$
satisfying $\frac{R}{r}<\sqrt{2}$, $\alpha^{-1}+\beta^{-1}=1$.

{\bf Lemma 3.1.} {\it  $-\frac{1}{2}(1-\frac{\alpha}{2})(\frac{\alpha}{2\pi r^2})^{\frac{\alpha}{2-\alpha}}\leq Min \Psi\leq
-\frac{1}{2}(1-\frac{\alpha}{2})(\frac{\alpha}{2\pi R^2})^{\frac{\alpha}{2-\alpha}}$.}

{\bf Proof.} Since $R^{-\alpha} |x|^\alpha\leq H(x)\leq r^{-\alpha} |x|^\alpha$, then by taking Fenchel conjugates we get
\bea \frac{1}{\beta}\left(\frac{r^\alpha}{\alpha}\right)^{\beta-1}|y|^\beta\leq G(y)\leq \frac{1}{\beta}\left(
\frac{R^\alpha}{\alpha}\right)^{\beta-1}|y|^\beta.\lb{3.1}\eea
Hence, we have \bea &&\Psi_r(u)\leq \Psi(u)\leq\Psi_R(u),\nn\\\Psi_r(u)&:=&\int_0^{\frac{1}{2}}\left(\frac{1}{2}J u \cdot
\Pi_\kappa u+\frac{1}{\beta}\left(\frac{r^\alpha}{\alpha}\right)^{\beta-1}|u|^\beta\right)dt\nn\\
\Psi_R(u)&:=&\int_0^{\frac{1}{2}}\left(\frac{1}{2}J u \cdot
\Pi_\kappa u+\frac{1}{\beta}\left(\frac{R^\alpha}{\alpha}\right)^{\beta-1}|u|^\beta\right)dt.\nn\eea
Then $Min \Psi_r\leq Min\Psi\leq Min\Psi_R$.

By the proof of Lemma 2.5 and replacing $H(x)$ with $r^{-\alpha} |x|^\alpha$, we obtain
\bea Min \Psi_r=-\frac{1}{2}(1-\frac{\alpha}{2})|\bar{x}|^\alpha r^{-\alpha},\lb{3.2}\eea
where $\bar{x}$ satisfies $Min \Psi_r=\Psi_r(\bar{x})$ and is a solution of:\bea\left\{\matrix{ \dot{x}=J\frac{\alpha}{r^\alpha}\frac{x}{|x|^{2-\alpha}}, \cr
x(\frac{1}{2})=P x(0).~~~~\cr }\right.\lb{3.3}\eea
By Lemma 2.5, $\bar{x}$ has minimal period 1 when $\bar{x}$ corresponds to a P-symmetric closed characteristic
or $\bar{x}$ has minimal period $\frac{1}{2}$ when $\bar{x}$ corresponds to a P-fixed closed characteristic.
By solving (\ref{3.3}), $\bar{x}$ has the form $\bar{x}(t)=e^{J\omega t}\bar{x}(0)$.
Then by (\ref{3.3}), we have \bea \omega=\frac{\alpha}{r^\alpha}|x|^{\alpha-2}, \lb{3.4}\eea
and $\omega=2\pi$, or $\omega=4\pi$.

When $\omega=2\pi$, we have $\Psi_r(\bar{x})=-\frac{1}{2}(1-\frac{\alpha}{2})(\frac{\alpha}{2\pi r^2})^{\frac{\alpha}{2-\alpha}}$ by (\ref{3.2}), (\ref{3.4}).

When $\omega=4\pi$, we have $\Psi_r(\bar{x})=-\frac{1}{2}(1-\frac{\alpha}{2})(\frac{\alpha}{4\pi r^2})^{\frac{\alpha}{2-\alpha}}$ by (\ref{3.2}), (\ref{3.4}).

Note that $Min \Psi_r=\Psi_r(\bar{x})$, then we must have $\omega=2\pi$ and
$\Psi_r(\bar{x})=-\frac{1}{2}(1-\frac{\alpha}{2})(\frac{\alpha}{2\pi r^2})^{\frac{\alpha}{2-\alpha}}$. Hence our Lemma holds.\hfill\hb

Recall that the action of a closed characteristic
$(\tau,y)$ is defined by (cf. P.190 of \cite{Eke1})\[
A(\tau,y)=\frac{1}{2}\int_0^{\tau}{(Jy\cdot \dot{y})}dt.\] Note that
$A(\tau,y)$ is a geometric quantity depending only on how many times
one runs around the closed characteristic. In fact, we can compute it as follows
\bea A(\tau,y)&=&\frac{1}{2}\int_0^{\tau}{(Jy\cdot JN_\Sigma(y))}dt\nn\\&=&\frac{1}{2}\int_0^{\tau}{1}dt\nn\\&=&\frac{\tau}{2}.\lb{3.5}\eea
Here we used (\ref{1.1}) and the fact that $\mathit{N}_\Sigma(y)\cdot y=1$.

Now we state a lemma dues to C. Croke and A. Weinstein which is useful for our proofs:

{\bf Lemma 3.2.}(cf. Theorem V.1.4 of \cite{Eke1}) {\it Assume $\Sigma$ is a $C^1$ hypersurface bounding a convex compact set $U$.
Suppose there is a point $x_{0}\in {\bf R}^{2n}$ such that
\[0<r \leq |x-x_{0}|,~\forall~x \in \Sigma\] Then, if $(\tau, y)$ is a closed characteristic on $\Sigma$, we have $A(\tau,y)\geq \pi r^2$.}

Now we restrict $\Psi$ on the subspace $X\equiv\{u=0\diamond u_2 \in L^2((0,\frac{1}{2}),
{\bf R}^{2n})\mid  u_2\in L^2((0,\frac{1}{2}), {\bf
R}^{2\kappa}), u(\frac{1}{2})=Pu(0),
\int_0^{\frac{1}{2}}u_2(t)dt=0\}$ of $L^2_\kappa\left(0,
\frac{1}{2}\right)$. Then by Lemma 2.2, $\Psi$ is bounded
from below on $X$, let $\Psi(\bar{u})=Min\{ \Psi(u)\mid u\in X\}$, then $\bar{u}$ corresponds to
a P-fixed closed characteristic via Proposition 2.1 and Lemma 2.2, precisely, $\bar{u}=\dot{x}$, $x={(\frac{\tau}{\alpha})}^{-\frac{1}{2-\alpha}} y(\tau
t)$, where $(\frac{\tau}{2},y)$ is a P-fixed closed characteristic. Hence, by Lemma 2.5 and (\ref{3.5}), we obtain\bea \Psi(\bar{u})=
-\frac{1}{2}(1-\frac{\alpha}{2})\left(\frac{4A(\frac{\tau}{2}, y)}{\alpha}\right)^{-\frac{\alpha}{2-\alpha}}.\lb{3.6}\eea
By Lemma 3.2, we have \bea A(\frac{\tau}{2}, y)\geq \pi r^2.\lb{3.7}\eea
Then (\ref{3.6}) and {\ref{3.7}) imply $\Psi(\bar{u})\geq
-\frac{1}{2}(1-\frac{\alpha}{2})\left(\frac{4\pi r^2}{\alpha}\right)^{-\frac{\alpha}{2-\alpha}}$.
Thus \bea \Psi(u)\geq
-\frac{1}{2}(1-\frac{\alpha}{2})\left(\frac{4\pi r^2}{\alpha}\right)^{-\frac{\alpha}{2-\alpha}} \lb{3.8}\eea for every $u\in X$. Now as Lemma V.2.5 of
\cite{Eke1}, we can
prove the following:

{\bf Lemma 3.3.} {\it When $\gamma\leq -\frac{1}{2}(1-\frac{\alpha}{2})\left(\frac{4\pi r^2}{\alpha}\right)^{-\frac{\alpha}{2-\alpha}}$, the $S^1$-action
is free in  $\Psi^\gamma$.}

{\bf Proof.} Firstly, by (\ref{3.8}), we have $\Psi^\gamma \cap X=\emptyset$. Let $u\in L^2_\kappa\left(0,
\frac{1}{2}\right)\setminus X$ has minimal period $\frac{1}{k}$ for some integer $k\geq2$. Note that by (\ref{2.4}), $k$ must be
odd, thus $k\geq3$. We have to prove that \bea \Psi(u)\geq
-\frac{1}{2}(1-\frac{\alpha}{2})\left(\frac{4\pi r^2}{\alpha}\right)^{-\frac{\alpha}{2-\alpha}}.\lb{3.9}\eea
In fact, define $v$ by \bea v(t):=k^{\frac{\alpha-1}{2-\alpha}}u(\frac{t}{k}).\nn\eea
Then $v\in L^2_\kappa\left(0,
\frac{1}{2}\right)$. Computing $\Psi(v)$, and using the $\beta$-homogeneous (since $H$ is $\alpha$-homogeneous), we get
\bea \Psi(v)=k^{\frac{\alpha}{2-\alpha}}\Psi(u).\nn\eea
Since $\Psi(v)\geq Min \Psi$, then by Lemma 3.1, we have\bea k^{\frac{\alpha}{2-\alpha}}\Psi(u)\geq
-\frac{1}{2}(1-\frac{\alpha}{2})(\frac{\alpha}{2\pi r^2})^{\frac{\alpha}{2-\alpha}}.\nn\eea
Then \bea \Psi(u)\geq
-\frac{1}{2}(1-\frac{\alpha}{2})(\frac{2k\pi r^2}{\alpha})^{-\frac{\alpha}{2-\alpha}}\geq
-\frac{1}{2}(1-\frac{\alpha}{2})\left(\frac{4\pi r^2}{\alpha}\right)^{-\frac{\alpha}{2-\alpha}}.\nn\eea
We complete our proof.\hfill\hb

In the following, we shall take:\bea \gamma>
-\frac{1}{2}(1-\frac{\alpha}{2})(\frac{\alpha}{2\pi R^2})^{\frac{\alpha}{2-\alpha}}.\lb{3.10}\eea
By Lemma 3.1, we have $\gamma>Min \Psi$, so $\Psi^\gamma$ is not empty. Since $2r^2>R^2$, we may also assume that\bea
 \gamma<
-\frac{1}{2}(1-\frac{\alpha}{2})(\frac{\alpha}{4\pi r^2})^{\frac{\alpha}{2-\alpha}}\lb{3.11}\eea
so the $S^1$-action on $\Psi^\gamma$ is free by Lemma 3.3. At this point we prove $\Psi^\gamma$
contains a copy of the Euclidian sphere $S^{2n-2\kappa-1}=\{x\in\R^{2n-2\kappa}\mid |x|^2=1\}$ with the standard
$S^1$-action. The standard
$S^1$-action on $S^{2n-2\kappa-1}$ is defined by:\bea \theta*x=e^{2\pi J_{2n-2\kappa}\theta}x.\nn\eea
We call a map $f:\R^{2n-2\kappa}\rightarrow L^2_\kappa\left(0,
\frac{1}{2}\right)$ equivariant if \bea f(e^{2\pi J_{2n-2\kappa}\theta}x)=\theta*f(x).\nn\eea

{\bf Lemma 3.4.} {\it There is an equivariant isometry $f:\R^{2n-2\kappa}\rightarrow L^2_\kappa\left(0,
\frac{1}{2}\right)$ and a number $\rho>0$ such that $f(\rho S^{2n-2\kappa-1})\subset\Psi^\gamma$.}

{\bf Proof.} We define $f(\xi)=u_\xi$ by $u_\xi(t)=(\sqrt{2}e^{2\pi J_{2n-2\kappa} t}\xi)\diamond 0$, $\forall \xi\in\R^{2n-2\kappa}$.
Then \bea \|f(\xi)\|=\left(\int_0^{\frac{1}{2}}|\sqrt{2}e^{2\pi J_{2n-2\kappa} t}\xi|^2 dt\right)^{\frac{1}{2}}=|\xi|,\nn\eea
that is, $f$ is an isometry. Since \bea f(e^{2\pi J_{2n-2\kappa} \theta}\xi)(t)=(\sqrt{2}e^{2\pi J_{2n-2\kappa} t}e^{2\pi J_{2n-2\kappa} \theta}\xi)\diamond 0
=(\sqrt{2}e^{2\pi J_{2n-2\kappa} (t+\theta)}\xi)\diamond 0=\theta*f(\xi),\nn\eea then $f$ is equivariant.

By direct computations, we have  $\Pi_\kappa
u_\xi(t)=(-\frac{\sqrt{2}}{2\pi}J_{2n-2\kappa}e^{2\pi
J_{2n-2\kappa}t}\xi)\diamond 0$, then \bea \frac{1}{2}\langle Ju_\xi, \Pi_\kappa u_\xi\rangle=-\frac{1}{4\pi}|\xi|^2.\nn\eea
and using (\ref{3.1}), we have \bea \Psi(u_\xi)&=&\int_0^{\frac{1}{2}}\left(\frac{1}{2}J u_\xi \cdot
\Pi_\kappa u_\xi+G(-Ju_\xi)\right)dt\nn\\&=&-\frac{1}{4\pi}|\xi|^2+\int_0^{\frac{1}{2}}G(-J(\sqrt{2}e^{2\pi J_{2n-2\kappa} t}\xi\diamond 0))dt
\nn\\&\leq& -\frac{1}{4\pi}|\xi|^2+\frac{1}{2}\cdot\frac{1}{\beta}\left(
\frac{R^\alpha}{\alpha}\right)^{\beta-1}|\sqrt{2}\xi|^\beta\nn\\&=&
\frac{1}{2}\left(-\frac{1}{4\pi}|\sqrt{2}\xi|^2+\frac{1}{\beta}\left(
\frac{R^\alpha}{\alpha}\right)^{\beta-1}|\sqrt{2}\xi|^\beta\right).\lb{3.12}\eea
The right-hand side of (\ref{3.12}) achieves its minimum when $-\frac{1}{2\pi}|\sqrt{2}\xi|+\left(
\frac{R^\alpha}{\alpha}\right)^{\beta-1}|\sqrt{2}\xi|^{\beta-1}=0$, which yields:\bea |\sqrt{2}\xi|=(\frac{1}{2\pi})^{\frac{\alpha-1}{2-\alpha}}
(\frac{\alpha}{ R^\alpha})^{\frac{1}{2-\alpha}}.\lb{3.13}\eea
Hence, by (\ref{3.10}) we have \be  \Psi(u_\xi)\leq
-\frac{1}{2}(1-\frac{\alpha}{2})(\frac{\alpha}{2\pi R^2})^{\frac{\alpha}{2-\alpha}}<\gamma,\lb{3.14}\ee
when (\ref{3.13}) holds. Let $\rho=\frac{1}{\sqrt{2}}(\frac{1}{2\pi})^{\frac{\alpha-1}{2-\alpha}}
(\frac{\alpha}{ R^\alpha})^{\frac{1}{2-\alpha}}$, we complete the proof of the lemma.

Recall that for a principal $U(1)$-bundle $E\rightarrow B$, the
Fadell-Rabinowitz index (cf. \cite{FaR1}) of E is defined to be
$\sup\{k \mid c_{1}(E)^{k-1}\neq 0\}$, where $c_1(E)\in H^2(B,{\bf
Q})$ is the first rational Chern class. For a $U(1)$-space, i.e., a
topological space $X$ with a $U(1)$-action, the Fadell-Rabinowitz
index $\hat{I}(X)$ is defined to be the index of the bundle $X\times
S^{\infty}\rightarrow X\times_{U(1)} S^{\infty}$, where
$S^{\infty}\rightarrow CP^{\infty}$ is the universal $U(1)$-bundle.

It is well known that $\hat{I}(S^{2n-2\kappa-1})=n-\kappa$, then by monotonicity of the Fadell-Rabinowitz index $\hat{I}$ and Lemma 3.4,
we have \bea \hat{I}(\Psi^\gamma)\geq n-\kappa.\lb{3.15}\eea
Now we define \begin{eqnarray}c_i = \inf\{\delta\in {\bf R}
\mid \hat{I}(\Psi^{\delta})\geq i\}.\nn\end{eqnarray}
Then \bea Min \Psi=c_1\leq c_2\leq\cdots\leq c_{n-\kappa}\leq\gamma<0.\lb{3.16}\eea
As Proposition V.3.3 in P.218 of \cite{Eke1}, we have

{\bf Proposition 3.5.} {\it Every $c_{i}$ is a critical value of
$\Psi$. If $c_{i} = c_{j}$ for some $i < j\leq n-\kappa$, then there are
infinitely many $S^1$-orbits of critical points on the level $c_i=c_j$.}

{\bf Proof.} For the reader's convenience, we sketch a brief proof here and refer to Sections V.2 and V.3
of \cite{Eke1} for related details.

By the proof of Theorem V.2.9 of \cite{Eke1}, if we replace $L_o^\beta$ and $\psi$ by $L_\kappa^2(0,\frac{1}{2})$
and $\Psi$ respectively, the Theorem V.2.9 of \cite{Eke1} also works.
Since the Fadell-Rabinowitz index  $\hat{I}$ has the properties of monotonicity, subadditivity, continuity which
are the only three properties of $I$ used in the proof of Proposition V.2.10 of \cite{Eke1}, then the
proof carries over verbatim of that of Proposition V.2.10 of \cite{Eke1}.

Now we prove the main theorems.

{\bf Proof of Theorem 1.1.} By (\ref{3.16}) and Proposition 3.5,
we find at least $n-\kappa$ distinct $S^1$-orbits of critical points $u_1, \cdots, u_{n-\kappa}$ of $\Psi$ with
$\Psi(u_i)<\gamma$. Using Lemma 3.1 and taking the infimum of $\gamma$ in formula (\ref{3.10}), we obtain:
\bea -\frac{1}{2}(1-\frac{\alpha}{2})(\frac{\alpha}{2\pi r^2})^{\frac{\alpha}{2-\alpha}}\leq \Psi(u_i)\leq
-\frac{1}{2}(1-\frac{\alpha}{2})(\frac{\alpha}{2\pi R^2})^{\frac{\alpha}{2-\alpha}}.\lb{3.17}\eea
By (\ref{3.8}) and Lemma 3.3, we know that every $u_i$ must corresponds to a prime $P$-symmetric closed characteristic,
denote it by $(\tau_i, y_i)$. Moreover by Lemma 2.5 and (\ref{3.5}), (\ref{3.17}), we have\bea \pi r^2\leq A(\tau_i, y_i)\leq \pi R^2.\lb{3.18}\eea
The proof is complete.

{\bf Proof of Theorem 1.2.} By the main result of \cite{EkL1} (cf. also Theorem V.2.1 of \cite{EkL1}),
there exist at least $\kappa$ geometrically distinct closed
characteristics on $\Sigma(\kappa)$, where $\Sigma(\kappa)=\{z\in {\bf R}^{2\kappa}\mid 0\diamond
z\in\Sigma\}$, 0 is the origin in ${\bf R}^{2n-2\kappa}$, then we obtain at least $\kappa$ geometrically distinct $P$-fixed closed
characteristics $\{(\tau_i, y_i)\mid n-\kappa+1\leq i\leq n\}$ on $\Sigma$ satisfying $\pi r^2\leq A(\tau_i, y_i)\leq \pi R^2$. Together
with Theorem 1.1, it proves Theorem 1.2.

{\bf Acknowledgements.} The authors would like to sincerely thank the referee for his/her
careful reading of the manuscript and valuable comments.

\bibliographystyle{abbrv}

\end{document}